\documentclass[11pt]{amsart}
\usepackage{latexsym,amssymb,amsmath}
\usepackage[all,ps,cmtip]{xy} 
\newtheorem{theo}{Theorem}
\newtheorem{coro}{Corollary}
\newtheorem{prop}{Proposition}
\newtheorem{lemm}{Lemma}

\begin{document}

\def\ot{\otimes}
\def\we{\wedge}
\def\wec{\wedge\cdots\wedge}
\def\op{\oplus}
\def\ra{\rightarrow}
\def\lra{\longrightarrow}
\def\fso{\mathfrak so}
\def\cO{\mathcal{O}}
\def\cS{\mathcal{S}}
\def\cL{\mathcal{L}}
\def\fsl{\mathfrak sl}
\def\fg{\mathfrak g}\def\fp{\mathfrak p}
\def\PP{\mathbb P}\def\QQ{\mathbb Q}\def\ZZ{\mathbb Z}\def\CC{\mathbb C}
\def\RR{\mathbb R}\def\HH{\mathbb H}\def\OO{\mathbb O}
\def\smc{\cdots }
\title[The canonical strip phenomenon]{The canonical strip phenomenon \\
for complete intersections in homogeneous spaces}

\author{L. Manivel}
\date{March 31, 2009}

\begin{abstract}
We show that a refined version of Golyshev's canonical strip hypothesis
does hold for the Hilbert polynomials of 
complete intersections in rational homogeneous spaces.  
\end{abstract}

\maketitle

\section{Introduction}
Golyshev \cite{golyshev} recently made some intriguing observations about 
Hilbert 
polynomials of canonically or anti-canonically polarized complex 
projective varieties. 
For a Fano variety $X$, the polynomial defined by $H_{-K_X}(k)=\chi(X,-kK_X)$
for $k\in\ZZ$ has the symmetry 
$$H_{K_X}(z)=(-1)^{\dim X}H_{K_X}(-z-1).$$ In particular 
its zeroes are symmetric with respect to the vertical line $Re(z)=-
\frac{1}{2}$. 
Golyshev suggested to consider the following series of hypothesis:
\medskip

\begin{tabular}{ll}
(CS) & The roots of $H_{-K_X}(z)$ belong to the {\it canonical strip}
$-1<Re(z)<0$. \\
(NCS) & The roots of $H_{-K_X}(z)$ belong to the {\it narrow 
canonical strip}. \\ 
(CL) &  The roots of $H_{-K_X}(z)$ belong to the {\it canonical line}
$Re(z)=-\frac{1}{2}$.
\end{tabular}

\medskip
The above-mentionned 
 {\it narrow canonical strip} is defined by the conditions that 
$$ -1 +\frac{1}{\dim X+1} < Re(z) <-\frac{1}{\dim X+1}.$$
We introduce the following variant for a Fano variety $X$ of index $\iota_X$. 
We define the {\it tight canonical strip} by the conditions 
$$-1+\frac{1}{i_X}<Re(z)<-\frac{1}{i_X},$$
and the corresponding hypothesis
\medskip

\begin{tabular}{cc}
(TCS) & The roots of $H_{-K_X}(z)$ belong to the {\it tight canonical strip}. 
\end{tabular}

\medskip
Golyshev noticed that (NCS) holds for Fano threefolds and that (CS) holds for
minimal threefolds of general type, as a consequence of well-known estimates
of their characteristic numbers. Shramov checked that (CL) holds for smooth 
Fano toric manifolds up to dimension four \cite{shr}. 

\smallskip
Note that we can express $H_{-K_X}(z)$ as follows. 
Consider the entire function 
$$\varphi(t)=\frac{t/2}{sh(t/2)}
=\sum_{t=0}^{\infty}B_n(\frac{1}{2})\frac{t^n}{n!},$$
where $B_n(u)$ denotes the $n$-th Bernouilli polynomial. This function 
defines a caracteristic class $\varphi(X)\in H^*(X,\QQ)$, and we can let 
$$\psi_k(X)=\int_X \varphi_k(X)\cup c_1(X)^{\dim X-k}.$$
Note that this is zero if $k$ is odd. It is then a direct consequence of the 
Grothendieck-Riemann-Roch theorem that 
$$H_{-K_X}(z-\frac{1}{2})=\sum_{k=0}^{\dim X}\psi_{\dim X-k}(X)
\frac{z^k}{k!}.$$
Write $\dim X=2m+\epsilon$, with $\epsilon\in\{0,1\}$. Then 
$H_{-K_X}(z-\frac{1}{2})=z^{\epsilon}P_X(z^2)$ for some 
polynomial of degree $m$, and the 
canonical line hypothesis (CL) can be translated into the 
hypothesis that the polynomial $P_X$ has only real and non-positive
roots.  

In this note we prove that:
\begin{enumerate}
\item (TCS) holds for rational homogeneous spaces of Picard number one;
\item (TCS) holds for all Fano complete intersections in rational 
homogeneous spaces of Picard number one;
\item (CL) holds for all general type complete intersections 
in rational homogeneous spaces of Picard number one.
\end{enumerate}

The first statement is a relatively straightforward application of the 
Weyl dimension formula and of the combinatorics of root systems. The next 
two statements can be deduced from the first one by induction.

\section{The tight canonical strip for rational homogeneous spaces}

Let $X=G/P$ be a complex rational homogeneous space with Picard number one. Here
$G$ denotes a simple affine algebraic group and $P$ a maximal parabolic
subgroup. Once we choose a maximal torus $T$ in $P$, we get a root system
$\Phi$ and its decomposition into positive and negative roots, 
$\Phi=\Phi_+\cup (-\Phi_+)$. Our (slightly unusual) convention will be
that negative roots are roots of $P$. Recall that the choice of $P$ 
(up to conjugation) is equivalent to the choice of a fundamental weight 
$\omega_0$, or equivalently, of a simple root $\alpha_0$.

Let $L$ be the ample
line bundle generating $Pic(X)$. It can be defined as the line bundle 
$L_{\omega_0}$ associated to the  fundamental weight $\omega_0$, considered 
as a character of $P$. The index $\iota_X$ of $X$ is defined by the 
identity $-K_X=\iota_X L$. 

By the Bott-Borel-Weil theorem, we know that 
$\Gamma(X,L^k)=V_{k\omega_0}$
is the irreducible $G$-module of highest weight $k\omega_0$, and that 
the higher cohomology groups vanish. We can then use Weyl's dimension
formula (see e.g. \cite{serre}) to express the Hilbert polynomial of $L$ as 
$$H_L(z)=\prod_{\alpha\in\Phi_+}\frac{(z\omega_0+\rho,\alpha)}{
(\rho,\alpha)}.$$
Here we have used an invariant pairing $(\; ,\; )$ on the weight lattice.
We will use the normalization defined  by the condition 
that $(\omega_0,\alpha_0)=1$, or equivalently, that $(\alpha_0,\alpha_0)=2$. 

We can decompose $H_L$ as the product of the 
polynomials $H_L^\ell$ defined as
 $$H_L^\ell(z)=\prod_{\substack{\alpha\in\Phi_+,\\ (\omega_0,\alpha)=\ell}}
\frac{\ell z+(\rho,\alpha)}{(\rho,\alpha)}=
\prod_j\Bigl(\frac{\ell z+j}{j}\Bigr)^{h_{\ell,j}},$$
where $h_{\ell,j}$ denotes the number of positive roots $\alpha\in\Phi_+$
such that $(\omega_0,\alpha)=\ell$ and $(\rho,\alpha)=j$. 

\medskip\noindent {\it Remark}. 
It is usually more natural to write down Weyl's dimension formula
in terms of the  coroots $\alpha^\vee$ attached to the roots $\alpha$.
Then $(\rho,\alpha^\vee)$ is an integer called the height of 
$\alpha$ and sometimes denoted $ht(\alpha)$. Since $\rho$ is equal 
to the sum of the fundamental weights, whose dual basis is that of 
the simple coroots, 
$ht(\alpha)$ can be computed as the sum of the 
coefficients of $\alpha$ over the basis of simple coroots.
See for example \cite{gw}. 

For our purposes, it seems necessary to work with roots rather than 
coroots. In the simply laced case this makes no difference, but in 
the non simply laced case we need to be careful with the fact that 
$(\rho,\alpha)$ can take non integer values. To avoid notational
complications, we will suppose in the sequel that $\fg$ is simply 
laced, but our results also hold in the non simply laced case, with 
essentially the same proofs.

\begin{prop} 
For each $\ell$, the sequence $h_{\ell,j}$ is symmetric and unimodal. 
That is 
\begin{eqnarray*}
(S)\qquad h_{\ell,j} &= &h_{\ell,\ell\iota_X-j}, \\
(U)\qquad h_{\ell,j} &\le &h_{\ell,j+1} \quad \mathit{if}\; 2j+1\le \ell\iota_X.
\end{eqnarray*}
\end{prop}

To prove this statement, we observe that there is a close connection between
the decomposition of $H_L$ into the product of the $H_L^\ell$'s, and the 
$\ZZ$-grading on $\fg=Lie(G)$ induced by $\omega_0$:
\begin{eqnarray*}
\fg  =\fg_{-\ell_{max}}\oplus\cdots\oplus \fg_{-i}\oplus
\cdots\oplus \fg_0\oplus\cdots\oplus \fg_i\oplus
\cdots\oplus\fg_{\ell_{max}}.
\end{eqnarray*}
Here $\fg_i$ denotes the sum of the root spaces associated to the 
roots having coefficient $i$ on $\alpha_0$ (plus the
Cartan subalgebra for $i=0$).
So the roots that contribute to $H_L^\ell$ are precisely those 
appearing in $\fg_\ell$, for $\ell\ge 1$. Now, we know 
(see e.g. \cite{rub})
that
\begin{itemize}
\item $\fg_0$ is a reductive subalgebra of $\fg$ of maximal rank, 
with rank one center;
\item the Dynkin diagram encoding the semi-simple part of $\fg_0$ can be 
deduced from the Dynkin diagram of $\fg$ by erasing the node corresponding
to the fundamental weight $\omega_0$, and the edges attached to this node;    
\item each $\fg_\ell$ is a simple $\fg_0$-module.
\end{itemize}
This implies that for each $\ell$ such that $\fg_\ell$ is non zero, 
there is a unique root $\gamma_\ell$ of $\fg$ which is a highest weight 
of $\fg_\ell$ considered as a $\fg_0$-module. Symmetrically, there is 
also a unique root $\beta_\ell$ of $\fg$ which is a lowest weight 
of $\fg_\ell$. Moreover, 
$$\beta_\ell=w_{00}(\gamma_\ell),$$
where $w_{00}$ denotes the longest element in the Weyl group of $\fg_0$
(not to be confused with the longest element in the Weyl group of $\fg$, 
usually denoted $w_0$), which is the subgroup of the Weyl group of 
$\fg$ generated by the simple reflections other than $s_{\alpha_0}$. 
Note in particular that 
$$(\rho,\beta_\ell+\gamma_\ell)=(\rho+w_{00}(\rho),\gamma_\ell).
$$
Note also that $\ell_{max}$ is the coefficient of the highest root $\psi$
of $\fg$ on $\alpha_0$. 

\begin{lemm}
We have $\rho+w_{00}(\rho)=\iota_X\omega_0$.
\end{lemm}

\proof Recall that $2\rho$ is the sum of the positive roots. 
The tangent bundle of $X=G/P$ is the homogeneous 
bundle defined by the $P$-module $\fg/\fp$, where $\fp=Lie(P)$. But
$$\fg/\fp \simeq \fg_1\oplus\cdots\oplus\fg_{\ell_{max}}.$$
This implies that the anticanonical bundle of $X$, being the determinant
of the tangent bundle, is defined by the weight $2\rho_X$ equal to the 
sum of the roots $\alpha$ such that $(\omega_0,\alpha) >0$. 

We can conclude that $2\rho-2\rho_X$ is the sum of the positive roots 
in $\fg_0$. In particular the action of $w_{00}$ takes it to its opposite. 
On the contrary, $2\rho_X=\iota_X\omega_0$ is not affected by the action of
$w_{00}$. Therefore
$$w_{00}(2\rho-2\rho_X)=w_{00}(2\rho)-2\rho_X=2\rho_X-2\rho,$$
and our claim follows. \qed

\medskip\noindent {\it Remark}. 
Recall that $2\rho$ is also the sum of the fundamental weights. 
For the root system of $\fg_0$, which is the set of roots of $\fg$ with 
zero coefficient on $\alpha_0$, the fundamental weights are of the form 
$\omega_i^0=\omega_i-a_i\omega_0$. Having coefficient zero on $\alpha_0$ 
is equivalent to being orthogonal to $\omega_0$, so 
$a_i=(\omega_i,\omega_0)/(\omega_0,\omega_0)$ and therefore
$2\rho_0=\sum_i\omega_i^0=2\rho-
\frac{(2\rho,\omega_0)}{(\omega_0,\omega_0)}\omega_0.$
Applying $w_{00}$, which maps $\rho_0$ to $-\rho_0$, we deduce that 
$\iota_X$ is also given by the simple formula
$$\iota_X=\frac{(2\rho,\omega_0)}{(\omega_0,\omega_0)}.$$
Observe that $(2\rho,\omega_0)=\sum_{\ell\ge 1}\ell\dim\fg_\ell$.

\medskip
We deduce from the previous lemma, and for each $\ell$, the identity   
$$(\rho,\beta_\ell+\gamma_\ell)=
(\iota_X\omega_0, \gamma_\ell)=\iota_X\ell.$$
Now, $w_{00}$ acts on the roots in $\fg_\ell$, sending $\beta_\ell$ to 
$\gamma_\ell$, and more  generally a root of height $k$ to a root of 
height $ht(\beta_\ell)+ht(\gamma_\ell)-k=\iota_X\ell-k$. This implies the
symmetry property (S). 

\smallskip
Finally, the unimodality property (U) is a special case of a general
property of weights of $\fg$-modules. Indeed, if $V_{\lambda}$
is the irreducible $\fg$-module of highest weight $\lambda$, then 
any weight of $V_{\lambda}$ is of the form $\lambda-\theta$ for some 
$\theta$ in the weight lattice, and the numbers of weights (counted 
with multiplicities) $\lambda-\theta$ such that $(\rho,\theta)$
is a given integer $\ell$, form a unimodular sequence. This 
general property gives (U) when applied to the irreducible and 
multiplicity-free $\fg_0$-module $\fg_\ell$ (see e.g. \cite{stanley}).\qed

\medskip
We can summarize our discussion by the following statement. Let 
$b_\ell$ denote the height of $\beta_{\ell}$, the unique smallest 
root such that $(\omega_0,\beta_{\ell})=\ell$. 

\begin{prop}
The Hilbert polynomial of the ample generator $L$ of the Picard group
of $X=G/P$ can be expressed as
$$H_L(z)=\prod_{\ell=1}^{\ell_{max}}\prod_{k=b_\ell}^{\ell\iota_X-b_\ell}
\Bigl(\frac{\ell z+k}{k}\Bigr)^{h_{\ell,k}},$$
where for each $\ell$ the sequence $h_{\ell,k}$ is symmetric and unimodal. 
\end{prop}

This implies that the (TCS) hypothesis holds for $X$. More precisely:

\begin{coro}
The zeroes of 
$H_{-K_X}(z)$ are contained in 
the real segment 
$$[-1+\frac{1}{\iota_X},-\frac{1}{\iota_X}].$$
\end{coro}

 
\medskip\noindent {\it Example}. Consider $X=E_6/P_4$ where, according 
to the notations of Bourbaki \cite{bou}, 
$P_4$ is the maximal parabolic subgroup 
associated with the simple root $\alpha_4$. This simple root corresponds
to the triple node of the Dynkin diagram of type $E_6$. In this case we 
have $\ell_{max}=3$ and the contributions of the three levels $\ell=1,2,3$ 
can be analyzed as follows.

\smallskip $\ell=1$: the extremal roots are
$$\beta_1=\begin{matrix} 0&0&1&0&0\\ & & 0 & &\end{matrix}, 
\qquad\gamma_1=\begin{matrix} 1&1&1&1&1\\ & & 1 & &\end{matrix}.$$
In particular $b_1=ht(\beta_1)=1$, and $ht(\beta_1)+ht(\gamma_1)=7=\iota_X$. 
There are eighteen roots contributing to $H_X^1$, and  
$$H_X^1(z)=\Bigl(\frac{z+1}{1}\Bigr)\Bigl(\frac{z+2}{2}\Bigr)^3
\Bigl(\frac{z+3}{3}\Bigr)^5\Bigl(\frac{z+4}{4}\Bigr)^5
\Bigl(\frac{z+5}{5}\Bigr)^3\Bigl(\frac{z+6}{6}\Bigr).$$

\smallskip $\ell=2$: the extremal roots are
$$\beta_2=\begin{matrix} 0&1&2&1&0\\ & & 1 & &\end{matrix}, 
\qquad\gamma_1=\begin{matrix} 1&2&2&2&1\\ & & 1 & &\end{matrix}.$$
In particular $b_2=ht(\beta_2)=5$, and $ht(\beta_2)+ht(\gamma_2)=14=2\iota_X$. 
There are nine roots contributing to $H_X^2$, and 
 $$H_X^2(z)=\Bigl(\frac{2z+5}{5}\Bigr)\Bigl(\frac{2z+6}{6}\Bigr)^2
\Bigl(\frac{2z+7}{7}\Bigr)^3\Bigl(\frac{2z+8}{8}\Bigr)^2
\Bigl(\frac{2z+9}{9}\Bigr).$$

\smallskip $\ell=3$: the extremal roots are
$$\beta_3=\begin{matrix} 1&2&3&2&1\\ & & 1 & &\end{matrix}, 
\qquad\gamma_3=\begin{matrix} 1&2&3&2&1\\ & & 2 & &\end{matrix}.$$
In particular $b_3=ht(\beta_3)=10$, and $ht(\beta_2)+ht(\gamma_2)=21=
3\iota_X$. 
These are the only two roots contributing to $H_X^3$, and 
 $$H_X^3(z)=\Bigl(\frac{3z+10}{10}\Bigr)\Bigl(\frac{3z+11}{11}\Bigr).$$

We conclude in particular that $X$ has dimension $29$ and degree 
$996584151214080$. 

\medskip\noindent {\it Remarks}. 

1. The simplest situation is when 
$\ell_{max}=1$, which occurs when $X$ is {\it cominuscule}. Then
$\beta_1=\alpha_0$, and $\gamma_1=\psi$ is the highest root. In particular 
we recover the known fact that $\iota_X=h+1$, where $h$ denotes the 
Coxeter number. Moreover we get the relation 
$$\dim X=(\omega_0,\omega_0)\iota_X.$$
The zeroes of the Hilbert polynomial in that case are just the 
$-j/\iota_X$, with $0<j<\iota_X$.

2. The Hilbert polynomials of the adjoint varieties (this is the special
case where $X=G/P\subset \PP(\fg)$ is the projectivization of a minimal
nilpotent orbit) have been investigated in \cite{lm1,lm3} from the 
perspective of Vogel and Deligne works on the {\it universal Lie algebra}. 
The related work \cite{lm2} explores some other cases related to certain 
series generalizing the lines of Freudenthal's magic square.

\section{Complete intersections}

In this section we prove that the (TCS) hypothesis holds for any complete 
intersection $Y$ in the rational homogeneous space $X=G/P$. More precisely, 
we will show that the Hilbert polynomial of $L_Y$ on $Y$ can be expressed as
$$H_{L_Y}(z)=H^0_{L_Y}(z)\prod_{\ell=1}^{\ell^Y_{max}}
\prod_{k=b_\ell(Y)}^{\ell\iota_Y-b_\ell(Y)}
\Bigl(\frac{\ell z+k}{k}\Bigr)^{h_{\ell,k}(Y)},$$
with the additional properties that:
\begin{enumerate}
\item for each $\ell$, the sequence $h{\ell,k}(Y)$ is symmetric and unimodal;
\item the zeroes of the polynomial $H^0_{L_Y}(z)$ all belong to
the line $Re(z)=-\frac{1}{2}$. 
\end{enumerate}

\smallskip
Note that when $\iota_Y\le 0$, the product in the right hand side of the 
previous identity is empty, since the integers $b_\ell(Y)$ will always 
be positive.

\medskip
We proceed by induction. Let $Z=Y\cap H_d$ be the transverse intersection 
of $Y$ with a hypersurface of degree $d$. Then $\iota_Z=\iota_Y-d$ and the 
Hilbert polynomial of $L_Z$ is simply 
$$H_{L_Z}(z)=H_{L_Y}(z)-H_{L_Y}(z-d).$$
For a given $\ell$, consider the polynomial
$$H^\ell_{L_Y}(z)=\prod_{k=b_\ell(Y)}^{\ell\iota_Y-b_\ell(Y)}
\Bigl(\frac{\ell z+k}{k}\Bigr)^{h_{\ell,k}(Y)}.$$
Then we can write $P^\ell_{L_Z}(z)=H^\ell_{L_Y}(z)-H^\ell_{L_Y}(z-d)$ as a
product $A^{\ell}(z)(B^{\ell}_+(z)-B^{\ell}_-(z))$, where
\begin{eqnarray*}
A^{\ell}(z) &= &\prod_{k=b_\ell(Y)}^{\ell\iota_Y-b_\ell(Y)-\ell d}
\Bigl(\frac{\ell z+k}{k}\Bigr)^{\min(h_{\ell,k}(Y),h_{\ell,k+\ell d}(Y))}, \\
B^{\ell}_+(z) &= &\prod_{k=b_\ell(Y)}^{\ell\iota_Y-b_\ell(Y)}
\Bigl(\frac{\ell z+k}{k}\Bigr)^{h_{\ell,k}(Y)
-\min(h_{\ell,k}(Y),h_{\ell,k+\ell d}(Y))}, \\
B^{\ell}_-(z) &= &\prod_{k=b_\ell(Y)}^{\ell\iota_Y-b_\ell(Y)}
\Bigl(\frac{\ell z+k}{k}\Bigr)^{h_{\ell,k+\ell d}(Y)
-\min(h_{\ell,k}(Y),h_{\ell,k+\ell d}(Y))}.
\end{eqnarray*}
Let $h_{\ell,k}(Z)=\min(h_{\ell,k}(Y),h_{\ell,k+\ell d}(Y))$. First observe 
that the symmetry of the sequence $h_{\ell,k}(Y)$ relative to the substitution
$k\mapsto \ell\iota_Y-k$ implies that of the sequence $h_{\ell,k}(Z)$
relative to $k\mapsto \ell\iota_Z-k$. Moreover, the sequence $h_{\ell,k}(Z)$
is again unimodal. Indeed, for $k\le \iota_Y\ell/2-d\ell/2=\iota_Z\ell/2$, 
\begin{eqnarray}
h_{\ell,k}(Z)=h_{\ell,k}(Y).
\end{eqnarray}
This is clear for $k\le \iota_Y\ell/2-d\ell$ since then 
$h_{\ell,k}(Y)\le h_{\ell,k+d\ell}(Y)$ by the unimodality for $Y$. 
For $\iota_Y\ell/2-d\ell\le k\le \iota_Y\ell/2-d\ell/2$, we observe that
by symmetry, 
$$h_{\ell,k+d\ell}(Y)=h_{\ell,\iota_Y\ell-k-d\ell}(Y)\ge h_{\ell,k}(Y)$$
since $k\le \iota_Y\ell-k-d\ell\le\iota_Y\ell/2$. We can therefore write
\begin{eqnarray*}
A^{\ell}(z) &= &\prod_{k=b_\ell(Z)}^{\ell\iota_Z-b_\ell(Z)}
\Bigl(\frac{\ell z+k}{k}\Bigr)^{h_{\ell,k}(Z)},
\end{eqnarray*}
where $b_\ell(Z)=b_\ell(Y)$ and the sequence $h_{\ell,k}(Z)$ is  
symmetric and unimodal. Beware that we can have $b_\ell(Z)>
\ell\iota_Z-b_\ell(Z)$, in which case $A^{\ell}(z)=1$ by convention. 

Now, using (1) we can rewrite $B^{\ell}_+$ and $B^{\ell}_-$ as
\begin{eqnarray*}
B^{\ell}_+(z) &= &\prod_{2k>\ell\iota_Z}
\Bigl(\frac{\ell z+k}{k}\Bigr)^{h_{\ell,k}(Y)-h_{\ell,k+\ell d}(Y)}, \\
B^{\ell}_-(z) &= &\prod_{2k<\ell\iota_Z}
\Bigl(\frac{\ell z+k}{k}\Bigr)^{h_{\ell,k+\ell d}(Y)-h_{\ell,k}(Y)}.
\end{eqnarray*}
Oberve that by symmetry $B^{\ell}_-(z)=\pm B^{\ell}_+(-\iota_Z-z)$. 
Moreover the zeroes of $B^{\ell}_+(z)$ (resp. $B^{\ell}_-(z)$) are located 
to the left (resp. to the right) of the line $Re(z)=-\iota_Z/2$. 
This implies that $|B^{\ell}_+(z)|>|B^{\ell}_-(z)|$ if $Re(z)>-\iota_Z/2$, 
and $|B^{\ell}_+(z)|<|B^{\ell}_-(z)|$ if $Re(z)<-\iota_Z/2$ (see the
proof of Lemma 2.2 in \cite{golyshev}). 

The same conclusion holds for the polynomial $Q(z)=H^0_{L_Y}(z)
-H^0_{L_Y}(z-d)$. Since the zeroes of $H^0_{L_Y}(z)$
are all supposed to lie on the line  $Re(z)=-\iota_Y/2$, we have 
$|H^0_{L_Y}(z)|>|H^0_{L_Y}(z-d)|$ if 
$Re(z)>-\iota_Z/2$, and $|H^0_{L_Y}(z)|
<|H^0_{L_Y}(z-d)|$ if $Re(z)<-\iota_Z/2$. 

Now we can conclude our analysis as follows: we can write
$$H_{L_Z}(z)=H^0_{L_Z}(z)\prod_{\ell=1}^{\ell^Z_{max}}
\prod_{k=b_\ell(Z)}^{\ell\iota_Z-b_\ell(Z)}
\Bigl(\frac{\ell z+k}{k}\Bigr)^{h_{\ell,k}(Z)},$$
where the polynomial $H^0_{L_Z}(z)$ is equal to 
$$H^0_{L_Y}(z)\prod_\ell B^{\ell}_+(z)-
H^0_{L_Y}(z-d)\prod_\ell B^{\ell}_-(z).$$
As we have seen, these two products can have the same moduli only if 
$Re(z)=-\iota_Z/2$, so the zeroes of $H^0_{L_Z}(z)$ have
to lie on that line.

Note that this analysis is still correct if $\iota_Z$ is negative, in 
which case we simply have $H_{L_Z}(z)=H^0_{L_Z}(z)$. We 
have proved:

\begin{theo}
Let $Y$ be a smooth complete intersection in a rational homogeneous space
$X=G/P$ with Picard number one. Then: 
\begin{enumerate}
\item If $Y$ has general type, all the zeroes of $H_{K_Y}(z)$ lie on the line 
$Re(z)=-\frac{1}{2}$, so that the (CL) hypothesis holds; 
\item If $Y$ is Fano, the zeroes of $H_{-K_Y}(z)$ lie either on the line 
$Re(z)=-\frac{1}{2}$ or on the real segment $[-1+\frac{1}{\iota_Y},
-\frac{1}{\iota_Y}]$. In particular, the (TCS) hypothesis holds. 
\item If $Y$ is Calabi-Yau, and $L_Y$ denotes the restriction to $Y$ of the 
ample generator of $Pic(G/P)$, then the zeroes of the Hilbert polynomial  
$H_{L_Y}(z)$ are purely imaginary. 
\end {enumerate}
\end{theo}
 
In the Fano case we can be more precise and describe explicitely the 
real zeros $H_{-K_Y}(z)$, which are all rational. Note that the (CL) hypothesis
does hold when the index $\iota_Y$ equals one or two. 

\section{Miscellanii}

\subsection{Branched coverings}
A variant of the previous remarks allows to draw 
similar conclusions for double coverings of rational homogeneous spaces,
as we did for their hypersurfaces. Let $X=G/P$ be as above, $L$ the ample 
generator of its Picard group, $H\subset X$ a general hypersurface in the 
linear system $|2dL|$. Let $\pi :Y\rightarrow X$ be a double-cover branched
over $H$. Then $K_Y=\pi^*(K_X+dL)$, so that $Y$ is Fano of index $\iota_Y
=\iota_X-d$ if $d<\iota_X$.  Since $\pi_*\cO_Y=\cO_X\oplus L^{-d}$, 
we can easily compute the Hilbert polynomial of $-K_Y$:
$$H_{-K_Y}(z)=H_{-K_X}(z)+H_{-K_X}(z-\frac{d}{\iota_Y}).$$
Up to the sign, this is the same formula as that giving the Hilbert 
polynomial of the hypersurface $H$. Since that sign did not matter in our
proof of Theorem 1 for hypersurfaces, we can conclude:

\begin{prop}
The zeroes of the Hilbert polynomial of $Y$ lie either on the line 
$Re(z)=-\frac{1}{2}$ or on the real segment $[-1+\frac{1}{\iota_Y}, 
-\frac{1}{\iota_Y}]$. In particular the (TCS) hypothesis holds for $Y$. 
\end{prop}

This extends to Calabi-Yau double covers (that is, $d=\iota_X$): 
the (CL) hypothesis does hold for the Hilbert polynomial of the 
polarization $\pi^*L$ of $Y$. 
 
\subsection{Complete intersections in abelian varieties}
We conclude this note by observing that the (CL) hypothesis is easily 
seen to hold for varieties of general type 
which are complete intersections in abelian varieties. 

Let $A$ be an abelian variety of dimension $n+c$, and $L_1,\ldots ,L_c$
be ample line bundles. Consider $X=H_1\cap\cdots\cap H_c\subset A$,
a transverse intersection of hypersurfaces $H_1\in |L_1|, \ldots ,
H_c\in |L_c|$. 

\begin{prop}
The canonical line hypothesis (CL) holds for $X$. 
\end{prop}

\proof The canonical bundle of $X$ is $K_X=(L_1+\cdots +L_c)|X$. 
Using the Koszul complex of $X$, we can compute its Hilbert polynomial as
\begin{eqnarray*}
H_{-K_X}(z) &= &\frac{1}{(n+c)!}\sum_\epsilon (-1)^{|\epsilon|}\bigl(
(z-\epsilon_1)L_1+\cdots +(z-\epsilon_c)L_c \bigr)^{n+c} \\
  &= &\sum_{\ell}\frac{(L_1^{\ell_1}\cdots L_c^{\ell_c})}{\ell_1!\cdots\ell_c!}
\sum_\epsilon (-1)^{|\epsilon|}(z-\epsilon_1)^{\ell_1}\cdots 
(z-\epsilon_c)^{\ell_c}.
\end{eqnarray*}
In these sums $\epsilon=(\epsilon_1,\ldots ,\epsilon_c)$ belongs to 
$\{0,1\}^c$,  $ |\epsilon|=\epsilon_1+\cdots +\epsilon_c$, and 
$\ell_1,\ldots ,\ell_c$ are non-negative integers of sum $n+c$. 
The polynomial $P_{\ell_1,\ldots ,\ell_c}(z)$ defined by the previous 
sum over $\epsilon$ factors as
$$P_{\ell_1,\ldots ,\ell_c}(z)=P_{\ell_1}(z)\cdots P_{\ell_c}(z),$$
where $P_\ell(z)=z^\ell-(z-1)^\ell$. There remains to observe that 
$P_\ell(z+\frac{1}{2})$ is a polynomial with non-negative coefficients,
of the same parity as $\ell$. Since the intersection coefficients 
$(L_1^{\ell_1}\cdots L_c^{\ell_c})$ are all positive, $H_{-K_X}(z)$
is therefore a polynomial in $z^2$ with positive coefficients, 
multiplied by $z$
if $n$ is odd. As we noticed in the introduction, 
this is enough to imply (CL). \qed

\bigskip

\providecommand{\bysame}{\leavevmode\hbox to3em{\hrulefill}\thinspace}
{\sc Laurent Manivel},

Institut Fourier, UMR 5582 (UJF-CNRS),

BP 74, 38402 St Martin d'H\`eres Cedex, France.

E-mail : {\tt Laurent.Manivel@ujf-grenoble.fr}

\end{document}